\documentclass[14pt]{article}
\usepackage{mathrsfs}
\usepackage{amsthm}
\usepackage{amsmath}
\usepackage{graphicx}
\usepackage{color}
\usepackage{amsfonts}
\usepackage{float}
\usepackage[text={140mm,210mm},left=45mm,vmarginratio=1:1]{geometry}
\newtheorem{theorem}{Theorem}[section]

\newtheorem{lemma}[theorem]{Lemma}

\newtheorem{corollary}[theorem]{Corollary}

\numberwithin{equation}{section}
\normalsize

\begin{document}
\title{\textbf{Phase transition for large dimensional contact process with random recovery rates on open clusters}}

\author{Xiaofeng Xue \thanks{\textbf{E-mail}: xuexiaofeng@ucas.ac.cn \textbf{Address}: School of Mathematical Sciences, University of Chinese Academy of Sciences, Beijing 100049, China.}\\ University of Chinese Academy of Sciences}

\date{}
\maketitle

\noindent {\bf Abstract:} In this paper we are concerned with
contact process with random recovery rates on open clusters of bond
percolation on $\mathbb{Z}^d$. Let $\xi$ be a positive random
variable, then we assigned i. i. d. copies of $\xi$ on the vertices
as the random recovery rates. Assuming that each edge is open with
probability $p$ and $\log d$ vertices are occupied at $t=0$, we
prove that the following phase transition occurs. When the infection
rate $\lambda<\lambda_c=1/({p{\rm E}\frac{1}{\xi}})$, then the
process dies out at time $O(\log d)$ with high probability as
$d\rightarrow+\infty$, while when $\lambda>\lambda_c$, the process
survives with high probability.


\section{Introduction}\label{section one}
In this paper, we are concerned with contact process with random
recover rates on open clusters of bond percolation in lattice. First
we introduce some notations. For each $d\geq 1$, we denote by
$\mathbb{Z}^d$ the $d$-dimensional lattice and denote by
$\mathbb{E}^d$ the set of edges on $\mathbb{Z}^d$. For any $x,y\in
\mathbb{Z}^d$, we denote by $x\sim y$ when there exits $e\in
\mathbb{E}^d$ connecting $x$ and $y$. Let $\{X(e)\}_{e\in
\mathbb{E}^d}$ be i. i. d. random variables such that
\begin{equation}\label{equ 1.1}
P(X(e)=1)=p=1-P(X(e)=0)
\end{equation}
for each $e\in \mathbb{E}^d$, where $p\in (0,1)$, then we denote by
$x\leftrightarrow y$ when and only when $x\sim y$ and the edge $e$
connecting $x$ and $y$ satisfies that $X(e)=1$. For each $x\in
\mathbb{Z}^d$, we define
\begin{equation}\label{equ 1.2}
\mathcal{N}(x)=\{y:y\leftrightarrow x\}
\end{equation}
as the set of neighbors of $x$. Please note that $\mathcal{N}(x)$ is
a random set depending on $\{X(e)\}_{e\in \mathbb{E}^d}$. Let $\xi$
be a random variable such that $P(\xi\geq 1)=1$ and
$\{\xi(x)\}_{x\in \mathbb{Z}^d}$ be i. i. d. random variables such
that $\xi(x)$ and $\xi$ have identical probability distributions. We
assume that $\{\xi(x)\}_{x\in \mathbb{Z}^d}$ and $\{X(e)\}_{e\in
\mathbb{E}^d}$ are independent.

The contact process is a Markov process with state space
\[
\mathcal{P}(\mathbb{Z}^d)=\{A:A\subseteq \mathbb{Z}^d\}.
\]
For any $t\geq 0$, we denote by $C_t$ the state of the process at
the moment $t$. After $\{X(e)\}_{e\in \mathbb{E}^d}$ and
$\{\xi(x)\}_{x\in \mathbb{Z}^d}$ are given, $\{C_t\}_{t\geq 0}$
evolves as follows.
\begin{equation}\label{equ 1.3 transition rate}
C_t\rightarrow
\begin{cases}
C_t\setminus \{x\} \text{~at rate~} \xi(x) & \text{\quad if~}x\in C_t,\\
C_t\cup \{x\} \text{~at rate~} \frac{\lambda}{2d}|\mathcal{N}(x)\cap
C_t| & \text{\quad if~}x\not\in C_t,
\end{cases}
\end{equation}
where $\lambda$ is a positive parameter called the infection rate.

The contact process describes the spread of an infection disease on
$\mathbb{Z}^d$. Vertices in $C_t$ are infected while vertices out of
$C_t$ are healthy. An infected vertex $x$ waits for an exponential
time with rate $\xi(x)$ to become healthy. An healthy vertex $y$ is
infected at rate proportional to the number of infected neighbors.

When $p=1$ and $P(\xi=1)=1$, the model turns into the classic
contact process introduced by Harris in \cite{Har1974}. The two
books \cite{Lig1985} and \cite{Lig1999} authored by Liggett give a
detailed survey of the study of contact process.

In recent years, the contact process on random graph generated from
the percolation model is a popular topic. In \cite{Ber2011},
Bertacchi, Lanchier and Zucca study the contact process on
$C_{\infty}\times K_N$, where $C_{\infty}$ is the infinite open
cluster of the site percolation and $K_N$ is the complete graph with
$N$ vertices. They give criterions to judge whether the process will
die out. In \cite{Chen2009}, Chen and Yao show that the complete
convergence theorem holds for contact process on open clusters of
$\mathbb{Z}^d\times \mathbb{Z}^+$. In \cite{Xue2016}, Xue shows that
the contact process on open clusters of oriented bond percolation in
$\mathbb{Z}^d$ has critical value approximately $1/(dp)$ as $d$
grows to infinity, where $p$ is the probability that an given edge
is open.

The study of contact process with random recovery rates dates back
to 1980s. In \cite{Bra1991}, Bramson, Durrett and Schonmann show
that the contact process with random recovery rates on
$\mathbb{Z}^1$ has an `intermediate phase' in which the process
survives but does not grow linearly. In \cite{Lig1992}, Liggett
studies contact process with random recovery rates and random
infection rates on $\mathbb{Z}^1$ and gives a sufficient condition
for the process to survive.

\section{Main results}\label{section two}
In this section we give main results of this paper. First we
introduce some definition and notations. For $d\geq 1$, we assume
that $\{X(e)\}_{e\in \mathbb{E}^d}$ and $\{\xi(x)\}_{x\in
\mathbb{Z}^d}$ are defined under a probability space
$(\Omega_d,\mathcal{F}_d,\mu_d)$. The expectation operator with
respect to $\mu_d$ is denoted by ${\rm E}_{\mu_d}$. For any
$\omega\in \Omega_d$ and $\lambda>0$, we denote by
$P_\lambda^\omega$ the probability measure of the contact process
$\{C_t\}_{t\geq 0}$ with infection rate $\lambda$ and recovery rates
$\{\xi(\omega,x)\}_{x\in \mathbb{Z}^d}$ on open clusters generated
from $\{X(\omega,e)\}_{e\in \mathbb{E}^d}$. $P_\lambda^\omega$ is
called the quenched measure. The expectation operator with respect
to $P_\lambda^\omega$ is denoted by ${\rm E}_\lambda^\omega$. For
each $d\geq 1$, we define
\[
P_{\lambda,d}(\cdot)=\int
P_\lambda^\omega(\cdot)~\mu_d(d\omega)={\rm
E}_{\mu_d}[P_\lambda^\omega(\cdot)],
\]
which is called the annealed measure. The expectation operator with
respect to $P_{\lambda,d}$ is denoted by ${\rm E}_{\lambda,d}$. When
there is no misunderstanding, we write $P_{\lambda,d}$ and ${\rm
E}_{\lambda,d}$ as $P_\lambda$ and $E_\lambda$. For any $A\subseteq
\mathbb{Z}^d$, we write $C_t$ as $C_t^A$ when $C_0=A$.

Now we can give our main result.

\begin{theorem}\label{theorem 2.1 main result}
For each $d\geq1$, let $A(d)$ be a subset of $\mathbb{Z}^d$ such
that $|A(d)|=\lceil\log d\rceil$, then for any
$\lambda<\lambda_c=1/(p{\rm E}\frac{1}{\xi})$, there exists
$c(\lambda)>0$ such that
\begin{equation}\label{equ 2.1}
\lim_{d\rightarrow+\infty}P_{\lambda,d}(C_{c(\lambda)\log
d}^{A(d)}\neq \emptyset)=0,
\end{equation}
while for any $\lambda>\lambda_c$,
\begin{equation}\label{equ 2.2}
\lim_{d\rightarrow+\infty}P_{\lambda,d}(\forall~t>0,C_t^{A(d)}\neq
\emptyset)=1.
\end{equation}
\end{theorem}

According to Theorem \ref{theorem 2.1 main result}, phase transition
occurs when the infection rate $\lambda$ grows from
$\lambda<\lambda_c$ to $\lambda>\lambda_c$, where
\[
\lambda_c=1/(p{\rm E}\frac{1}{\xi}).
\]

We denote by $O$ the origin of $\mathbb{Z}^d$. According to the
basic coupling of spin systems (See Section 3.1 of \cite{Lig1985}),
it is easy to see that for any $\lambda_1\geq\lambda_2$
\[
P_{\lambda_1}(\forall~t>0, C_t^O\neq \emptyset)\geq
P_{\lambda_2}(\forall~t>0, C_t^O\neq \emptyset).
\]
Therefore, for each $d\geq 1$, the following definition is
reasonable.
\begin{equation}\label{equ 2.3}
\lambda_d=\sup\{\lambda:P_{\lambda,d}(\forall~t>0, C_t^O\neq
\emptyset)=0\}.
\end{equation}
According to Theorem \ref{theorem 2.1 main result}, we have the
following corollary.
\begin{corollary}\label{corollary 2.2}
$\lambda_d$ is as that defined in \eqref{equ 2.3}, then
\[
\limsup_{d\rightarrow+\infty}\lambda_d\leq \lambda_c=1/(p{\rm
E}\frac{1}{\xi}).
\]
\end{corollary}

When $p=\xi=1$, Corollary \ref{corollary 2.2} shows that
$\limsup_{d\rightarrow+\infty}\lambda_d\leq 1$. A stronger
conclusion such that $\lim_{d\rightarrow+\infty}\lambda_d=1$ for the
classic contact process is shown by Holley and Liggett in
\cite{Hol1981}. In \cite{Xue2014}, Xue shows that the critical value
of high dimensional threshold one contact process has similar
asymptotic behavior.

The proof of Theorem \ref{theorem 2.1 main result} is divided into
two sections. In Section 3, we will prove Equation \eqref{equ 2.1}.
The proof is inspired by the approach of graphical representation
introduced by Harris in \cite{Har1978}. In Section 4, we will prove
Equation \eqref{equ 2.2}. The proof is inspired by the approach
introduced by Kesten in \cite{Kesten1990} to study the asymptotic
behavior of the critical probability of high dimensional percolation
model.

\section{Subcritical case} \label{section 3}
In this section we give the proof of \eqref{equ 2.1}. First we
introduce the graphical representation of the process
$\{C_t\}_{t\geq 0}$. We consider the graph $\mathbb{Z}^d\times
[0,+\infty)$. In other words, we erect a time arrow on each vertex
on $\mathbb{Z}^d$. After $\{\xi(x)\}_{x\in \mathbb{Z}^d}$ and
$\{X(e)\}_{e\in \mathbb{E}^d}$ are given, we assume that
$\{Y_x(t)\}_{t\geq 0}$ is a Poisson process with rate $\xi(x)$ and
$\{U_{(x,y)}(t)\}_{t\geq 0}$ is a Poisson process with rate
$\lambda/(2d)$ for each $x,y\in \mathbb{Z}^d,x\sim y$. We assume
that all these Poisson processes are independent. Please note that
we care about the order of $x$ and $y$, hence $U_{(x,y)}\neq
U_{(y,x)}$. For any event time $t$ of $Y_x$, we put a `$\Delta$' at
$(x,t)$. For any event time $s$ of $U_{(x,y)}$, we put an arrow
`$\rightarrow$' from $(x,s)$ to $(y,s)$. For $x,y\in \mathbb{Z}^d$
and $t>0$, we say that there is an infection path from $(x,0)$ to
$(y,t)$ when there exist $x=x_0\sim x_1\sim x_2\sim\ldots\sim x_n=y$
and $0=t_{-1}<t_0<t_1<\ldots<t_n=t$ satisfying all the following
three conditions.

(1) For $0\leq i\leq n-1$, there is an arrow from $(x_i,t_i)$ to
$(x_{i+1},t_i)$.

(2) For $0\leq i\leq n$, there is no `$\Delta$' on
$\{x_i\}\times(t_{i-1},t_i]$.

(3) For $0\leq i\leq n-1$, $X(x_i,x_{i+1})=1$.

Please note that we write $X(e)$ as $X(x,y)$ when $e$ connecting $x$
and $y$. The following figure gives an example of infection path.

\begin{figure}[H]
  \centering
  \includegraphics[height=0.5\textwidth=0.5]{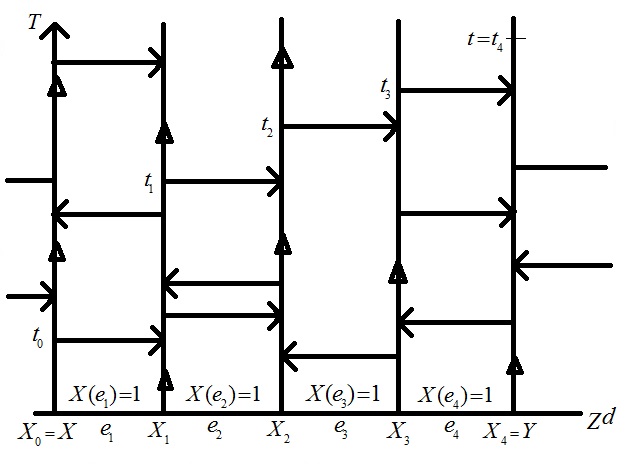}
  \caption{Infection path}
  \label{graph 3.1}
\end{figure}

In Figure \ref{graph 3.1}, there is an infection path from $(x_0,0)$
to $(x_4,t)$. According to the transition rates function of
$\{C_t\}_{t\geq 0}$ given by \eqref{equ 1.3 transition rate}, it is
easy to see that
\begin{equation}\label{equ 3.1}
C_t^A=\{y:\text{~for some~}x\in A,\text{~there is an infection path
from~}(x,0)\text{~to~}(y,t)\}
\end{equation}
for any $A\subseteq \mathbb{Z}^d$ in the sense of coupling. By
Equation \eqref{equ 3.1},
\[
C_t^A=\cup_{x\in A}C_t^x
\]
in the sense of coupling. Therefore, for any finite $A$,
\begin{equation}\label{equ 3.2}
P_\lambda^\omega(C_t^A\neq \emptyset)\leq \sum_{x\in
A}P_\lambda^\omega(C_t^x\neq \emptyset)
\end{equation}
and
\begin{equation}\label{equ 3.3}
P_{\lambda,d}(C_t^A\neq \emptyset)\leq |A|P_{\lambda,d}(C_t^O\neq
\emptyset).
\end{equation}

For later use, we divided the infection paths into several types.
For each $n\geq 0$, we define
\[
B_n=\{(x_0,x_1,x_2,\ldots,x_n):O=x_0\sim x_1\sim x_2\sim\ldots\sim
x_n\}
\]
as the set of path starting at $O$ with length $n$. For
$\vec{x}=(O,x_1,x_2,\ldots,x_n)\in B_n$ and $n$ positive integers
$j_0,j_1,\ldots,j_{n-1}$, we say that $\vec{x}$ is an infection path
with type $(j_0,j_1,\ldots,j_{n-1})$ at the moment $t$ when there
exists $0=t_{-1}<t_0<t_1<\ldots<t_{n}=t$ such that

(1) For $0\leq i\leq n-1$, there is an arrow from $(x_i,t_i)$ to
$(x_{i+1},t_i)$.

(2) For $0\leq i\leq n$, there is no `$\Delta$' on
$\{x_i\}\times(t_{i-1},t_i]$.

(3) For $0\leq i\leq n-1$, $X(x_i,x_{i+1})=1$.

(4) For $0\leq i\leq n-1$, $t_i$ is the $j_i$th event time of
$U_{(x_i,x_{i+1})}$ after the moment $t_{i-1}$.

The following figure gives an example.

\begin{figure}[H]
  \centering
  \includegraphics[height=0.5\textwidth=0.5]{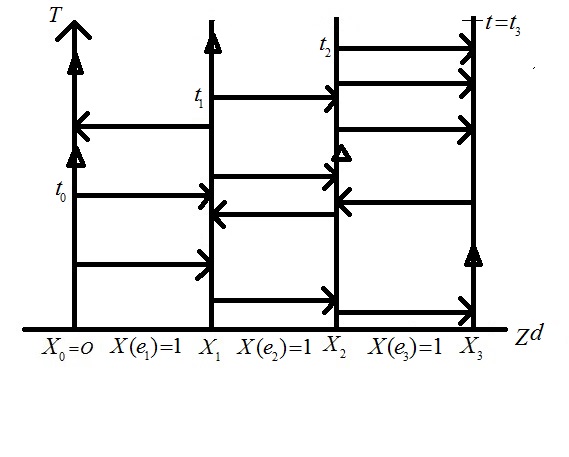}
  \caption{Infection path type}
  \label{graph 3.2}
\end{figure}

In Figure \ref{graph 3.2}, $(O,x_1,x_2,x_3)$ is an infection path
with type $(2,2,2)$ at the moment $t$. Please note that an infection
path may be with more than one type. In Figure \ref{graph 3.2},
$(O,x_1,x_2,x_3)$ is also with the type $(2,2,1)$. We use
$A(\vec{x},j_0,j_1,\ldots,j_{n-1},t)$ to denote the event that
$\vec{x}$ is an infection path with type $(j_0,j_1,\ldots,j_{n-1})$
at the moment $t$. We define
\[
A(\vec{x},t)=\bigcup_{j_0,j_1,\ldots,j_{n-1}}A(\vec{x},j_0,j_1,\ldots,j_{n-1},t)
\]
as the vent that $\vec{x}$ is an infection path at the moment $t$.
For later use, we define
\[
D_n=\{\vec{x}=(x_0,x_1,\ldots,x_n)\in B_n:x_i\neq x_j \text{~for
any~}0\leq i<j\leq n\}.
\]

After all the above prepared work, we give the proof of Equation
\eqref{equ 2.1}.

\proof[Proof of Equation \ref{equ 2.1}] We let
$c(\lambda)=\frac{1}{2\lambda}$ and $t=c(\lambda)\log d$. By
Equation \eqref{equ 3.1},
\begin{align}\label{equ 3.4}
P_{\lambda,d}(C_t^{O}\neq \emptyset)\leq&
\sum_{n=0}^{+\infty}\sum_{\vec{x}\in
D_n}\sum_{j_0,j_1,\ldots,j_{n-1}\geq
1}P_{\lambda,d}\big(A(\vec{x},j_0,\ldots,j_{n-1},t)\big)\\
&+\sum_{n=0}^{+\infty}\sum_{\vec{x}\in B_n\setminus
D_n}P_{\lambda,d}\big(A(\vec{x},t)\big).\notag
\end{align}
First we deal with $\sum_{n=0}^{+\infty}\sum_{\vec{x}\in
B_n\setminus D_n}P_{\lambda,d}\big(A(\vec{x},t)\big)$ .It is easy to
see that $B_n=D_n$ for $n=0,1$, hence
\begin{equation}\label{equ 3.5}
\sum_{n=0}^{+\infty}\sum_{\vec{x}\in B_n\setminus
D_n}P_{\lambda,d}\big(A(\vec{x},t)\big)=\sum_{n=2}^{+\infty}\sum_{\vec{x}\in
B_n\setminus D_n}P_{\lambda,d}\big(A(\vec{x},t)\big).
\end{equation}
For each $n\geq 2$, it is easy to see that
\begin{equation}\label{equ 3.6}
|B_n\setminus D_n|\leq {n+1 \choose 2}(2d)^{n-1}.
\end{equation}
This is because for any $\vec{x}\in B_n\setminus D_n$, there exists
$i\neq j$ such that $x_i=x_j$ and for each $l\neq 0,i,j$, $x_l$ has
at most $2d$ choices. Let $\{N(s)\}_{s\geq 0}$ be a Poison process
with rate $\frac{\lambda}{2d}$, then we claim that
\begin{equation}\label{equ 3.7}
P_{\lambda,d}(A(\vec{x},t))\leq P(N(t)\geq n)
\end{equation}
for each $\vec{x}\in B_n\setminus D_n$. To explain \eqref{equ 3.7},
we denote $T_0$ the first event time of $U_{(x_0,x_1)}$ and $T_l$
the first event time of $U_{(x_l,x_{l+1})}$ after the moment
$T_{l-1}$ for $1\leq l\leq n-1$. Therefore, $\{T_j-T_{j-1}\}_{0\leq
j\leq n-1}$ are i. i. d. exponential times with rate
$\frac{\lambda}{2d}$, where $T_{-1}=0$. As a result,
\begin{equation}\label{equ 3.8}
P(T_{n-1}<t)=P(N(t)\geq n).
\end{equation}
If $\vec{x}$ is an infection path at the moment $t$, then there
exists $t_0<t_1<\ldots<t_{n-1}<t_n=t$ such that $t_i$ is an event
time of $U_{(x_i,x_{i+1})}$ at which $x_i$ infects $x_{i+1}$.
According to the definition of $\{T_j\}_{0\leq j\leq n-1}$,
\[
T_j\leq t_j
\]
for each $1\leq j\leq n-1$. As a result,
\begin{equation}\label{equ 3.9}
P_{\lambda,d}(A(\vec{x},t))\leq P_{\lambda,d}(T_{n-1}<t).
\end{equation}
Equation \eqref{equ 3.7} follows from \eqref{equ 3.8} and \eqref{equ
3.9}.

For $\theta>0$,
\[
P(N(t)\geq n)\leq e^{-\theta n}{\rm E}e^{\theta
N(t)}=e^{\frac{\lambda t}{2d}(e^\theta-1)-\theta n}.
\]
For sufficiently large $d$, $n\geq 2$ and $t=\frac{1}{2\lambda}\log
d$, we choose $\theta=\log \frac{2dn}{\lambda t}$, then according to
the Stirling's Formula,
\begin{equation}\label{equ 3.10}
P(N(t)\geq n)\leq \frac{e^ne^{-\frac{\lambda t}{2d}}(\lambda
t)^n}{(2dn)^n}\leq \frac{M\sqrt{n}}{n!}\big(\frac{\lambda
t}{2d}\big)^n,
\end{equation}
where $M$ is a constant which does not depend on $\lambda,t,d$ and
$n$.

By Equation \eqref{equ 3.5}, \eqref{equ 3.6}, \eqref{equ 3.7} and
\eqref{equ 3.10}, when $t=\frac{1}{2\lambda}\log d$,
\begin{align}\label{equ 3.11}
\sum_{n=0}^{+\infty}\sum_{\vec{x}\in B_n\setminus
D_n}P_{\lambda,d}\big(A(\vec{x},t)\big)&\leq
M\sum_{n=2}^{+\infty}\frac{(n+1)n(2d)^{n-1}\sqrt{n}}{2n!}\big(\frac{\lambda
t}{2d}\big)^n\notag\\
&\leq \frac{M_1}{2d}\sum_{n=2}^{+\infty}\frac{n^3(\lambda
t)^n}{n!}\notag\\
&\leq \frac{M_2(\lambda t)^2+M_3(\lambda t)^3e^{\lambda t}}{2d}\notag\\
&=\frac{M_4\log^2 d}{d}+\frac{M_5 \log^3 d}{\sqrt{d}}
\end{align}
where $M_1, M_2, M_3, M_4, M_5$ are constant does not depend on
$\lambda,t,n$ and $d$.

\quad

Now we deal with
$P_{\lambda,d}(A(\vec{x},j_0,j_1,\ldots,j_{n-1},t))$ for $\vec{x}\in
D_n$. For given $\omega\in \Omega_d$, we let
$S_0,S_1,\ldots,S_{n-1}$ be independent exponential times with rates
$\xi(x_0,\omega),\xi(x_1,\omega),\ldots,\xi(x_n,\omega)$
respectively. We let $\{W_i\}_{1\leq i\leq +\infty}$ be i. i. d.
exponential times with rate $\frac{\lambda}{2d}$ and $\{V_m\}_{0\leq
m\leq n-1}$ be independent random variables such that $V_m$ and
$\sum_{l=1}^{j_m}W_l$ have identical probability distributions for
$0\leq m\leq n-1$. Then, according to the definition of
$A(\vec{x},j_0,j_1,\ldots,j_{n-1},t)$,
\begin{align}\label{equ 3.12}
&P_\lambda^\omega(A(\vec{x},j_0,j_1,\ldots,j_{n-1},t))\notag\\
&=p^nP(\sum_{m=0}^{n-1}V_m\leq t, S_m\geq V_m, 0\leq m\leq n-1,
S_n\geq t-\sum_{m=0}^{n-1}V_m).
\end{align}
Please note that the factor $p^n$ in Equation \eqref{equ 3.12} is
the probability that $X(x_i,x_{i+1})=1$ for $0\leq i\leq n-1$, since
$\{X(x_i,x_{i+1})\}_{0\leq i\leq n-1}$ are i. i. d. for
$\vec{x}=(x_0,x_1,\ldots,x_n)\in D_n\texttt{\texttt{}}$.

Since
\[
P(S_m\geq V_m|V_0,\ldots,V_{n-1})=e^{-\xi(x_m)V_m}
\]
and $V_m$ has probability density
\[
p(t_m,j_m)=\frac{t_m^{j_m-1}e^{-\frac{\lambda}{2d}t_m}}{(j_m-1)!}(\frac{\lambda}{2d})^{j_m}
\]
for $0\leq m\leq n-1$,
\begin{align}\label{equ 3.13}
&P_\lambda^\omega(A(\vec{x},j_0,j_1,\ldots,j_{n-1},t))\notag\\
&=p^n\int\limits_{\sum_{m=0}^{n-1}t_m\leq
t}\prod_{m=0}^{n-1}p(t_m,j_m)e^{-\sum\limits_{m=0}^{n-1}\xi(x_m)t_m}e^{-\xi(x_n)(t-\sum\limits_{m=0}^{n-1}t_m)}~dt_0~dt_1...~dt_m.
\end{align}
By Equation \eqref{equ 3.13} and direct calculation, it is not
difficult to check that
\begin{align}\label{equ 3.14}
&\sum\limits_{j_0,\ldots,j_{n-1}\geq
1}P_\lambda^\omega(A(\vec{x},j_0,j_1,\ldots,j_{n-1},t))\notag\\
&=(\frac{\lambda
p}{2d})^n\int\limits_{\sum\limits_{m=0}^{n-1}t_m\leq
t}e^{-\sum\limits_{m=0}^{n-1}\xi(x_m)t_m}e^{-\xi(x_n)(t-\sum\limits_{m=0}^{n-1}t_m)}~dt_0~dt_1\ldots~dt_{n-1}\notag\\
&=(\frac{\lambda
p}{2d})^n\frac{1}{\prod\limits_{m=0}^{n-1}\xi(x_m)}P(\sum\limits_{m=0}^{n-1}S_m\leq
t,S_n\geq t-\sum_{m=0}^{n-1}t_m).
\end{align}
The calculation is a little tedious, we omit the details. Let
$\{\alpha(t)\}_{t\geq 0}$ be a Poison process with rate $1$, then
\begin{equation}\label{equ 3.15}
P(\sum\limits_{m=0}^{n-1}S_m\leq t,S_n\geq
t-\sum_{m=0}^{n-1}S_m)\leq P(\sum\limits_{m=0}^nS_m\geq t)\leq
P(\alpha(t)\leq n+1),
\end{equation}
since $\xi(x)\geq 1$ for each $x\in \mathbb{Z}^d$. By Equation
\eqref{equ 3.14} and \eqref{equ 3.15},
\begin{align}\label{equ 3.16}
&\sum_{n=0}^{+\infty}\sum_{\vec{x}\in
D_n}\sum_{j_0,j_1,\ldots,j_{n-1}\geq
1}P_{\lambda,d}\big(A(\vec{x},j_0,\ldots,j_{n-1},t)\big)\notag\\
&\leq \sum_{n=0}^{+\infty}\sum_{\vec{x}\in D_n}(\frac{\lambda
p}{2d})^n{\rm
E}_{\mu_d}\Big[\frac{1}{\prod\limits_{m=0}^{n-1}\xi(x_m)}\Big]P(\alpha(t)\leq
n+1).
\end{align}
For $\vec{x}\in D_n$, $\{x_i\}_{0\leq i\leq n}$ are different with
each other, therefore
\begin{equation}\label{equ 3.17}
{\rm
E}_{\mu_d}\Big[\frac{1}{\prod\limits_{m=0}^{n-1}\xi(x_m)}\Big]=({\rm
E\frac{1}{\xi}})^n.
\end{equation}
Since each vertex on $\mathbb{Z}^d$ has $2d$ neighbors,
\begin{equation}\label{equ 3.18}
|D_n|\leq (2d)^n.
\end{equation}
By Equation \eqref{equ 3.16}, \eqref{equ 3.17} and \eqref{equ 3.18},
\begin{align}\label{equ 3.19}
\sum_{n=0}^{+\infty}\sum_{\vec{x}\in
D_n}\sum_{j_0,j_1,\ldots,j_{n-1}\geq
1}P_{\lambda,d}\big(A(\vec{x},j_0,\ldots,j_{n-1},t)\big)\leq
\sum_{n=0}^{+\infty}(\lambda p{\rm E}\frac{1}{\xi})^nP(\alpha(t)\leq
n+1).
\end{align}
For $n\leq \lfloor\frac{1}{4\lambda}\log d\rfloor-1$,
\begin{equation}\label{equ 3.20}
P(\alpha(t)\leq n+1)\leq P(\alpha(t)\leq t/2)\leq
2^{\frac{t}{2}}{\rm
E}(\frac{1}{2})^{\alpha(t)}=(\sqrt{\frac{2}{e}})^t.
\end{equation}
By Equation \eqref{equ 3.19} and \eqref{equ 3.20}, when
$\lambda<\lambda_c$,
\begin{align}\label{equ 3.21}
&\sum_{n=0}^{+\infty}\sum_{\vec{x}\in
D_n}\sum_{j_0,j_1,\ldots,j_{n-1}\geq
1}P_{\lambda,d}\big(A(\vec{x},j_0,\ldots,j_{n-1},t)\big)
\notag\\
&\leq \lfloor\frac{1}{4\lambda}\log d\rfloor
(\sqrt{\frac{2}{e}})^t+\sum\limits_{n=\lfloor\frac{1}{4\lambda}\log
d\rfloor}^{+\infty}(\lambda p{\rm
E}\frac{1}{\xi})^n\notag\\
&=\lfloor\frac{1}{4\lambda}\log d\rfloor
\big[(\frac{2}{e})^{\frac{1}{4\lambda}}\big]^{\log d}+\frac{(\lambda
p {\rm E}\frac{1}{\xi})^{\lfloor\frac{1}{4\lambda}\log
d\rfloor}}{1-\lambda p{\rm E}\frac{1}{\xi}}.
\end{align}

By Equation \eqref{equ 3.4}, \eqref{equ 3.11} and \eqref{equ 3.21},
when $\lambda<\lambda_c$,
\begin{equation}\label{equ 3.22}
\lim_{d\rightarrow+\infty}(\log
d)P_{\lambda,d}(C_{\frac{1}{2\lambda}\log d}^O\neq \emptyset)=0.
\end{equation}
Equation \eqref{equ 2.1} follows from Equation \eqref{equ 3.3} and
\eqref{equ 3.22} with $c(\lambda)=\frac{1}{2\lambda}$.

\qed

\section{Supcritical case} \label{section 4}
In this section we give the proof of Equation \eqref{equ 2.2}. Our
proof is inspired by a technique introduced in \cite{Kesten1990},
where $\lim_{d\rightarrow+\infty}2dp_c(d)=1$ is shown for the
critical probability $p_c(d)$ of $d$-dimensional site percolation.
First we introduce some notations. For each $d\geq 1$, we define
\[
N(d)=\lfloor\frac{\log d}{2\log\log d}\rfloor.
\]
We write $N(d)$ as $N$ when there is no misunderstanding. For $1\leq
i\leq d$, we define
\[
e_i=(0,~\ldots,~0,\mathop 1\limits_{i \text{th}},0,~\ldots,~0).
\]
For $x\in \mathbb{Z}^d$ and integer $k\geq 1$, we define
\begin{align}\label{equ 4.1 definition}
F_k(x)=&\Big\{\vec{x}=(x_0,x_1,\ldots,x_{kN-1}): x=x_0\sim
x_1\sim\ldots\sim x_{kN-1}, \text{~for~} N|(j+1), \notag\\
& x_{j+1}-x_j\in\{e_k:k>d-\lfloor\frac{d}{N}\rfloor\},
\text{~for~}N\not|(j+1),\\
&x_{j+1}-x_j\in \{\pm e_k: k\leq
d-\lfloor\frac{d}{N}\rfloor\}\Big\}\notag.
\end{align}
For each $\vec{x}\in F_k(x)$, we denote by $I_k(\vec{x})$ the event
that there exists $t>0$ such that $\vec{x}$ is an infection path
with type $(1,1,\ldots,1)$ at the moment $t$. On the event
$\bigcap\limits_{k=1}^{+\infty}\bigcup\limits_{x\in
A}\bigcup\limits_{\vec{x}\in F_k(x)}I_k(\vec{x})$, there exists
infection path starting at some vertex in $A$ and ending at vertex
with arbitrary large norm, hence the process will not die out.
Therefore,
\begin{align}\label{equ 4.2}
P_{\lambda,d}(C_t^A\neq \emptyset,\forall~ t>0)&\geq
P_{\lambda,d}\Big(\bigcap\limits_{k=1}^{+\infty}\bigcup\limits_{x\in
A}\bigcup\limits_{\vec{x}\in
F_k(x)}I_k(\vec{x})\Big)\notag\\
&=\lim\limits_{k\rightarrow+\infty}P_{\lambda,d}\Big(\bigcup\limits_{x\in
A}\bigcup\limits_{\vec{x}\in F_k(x)}I_k(\vec{x})\Big).
\end{align}
To deal with $P_{\lambda,d}\Big(\bigcup\limits_{x\in
A}\bigcup\limits_{\vec{x}\in F_k(x)}I_k(\vec{x})\Big)$ later, we
need the following lemma.
\begin{lemma}\label{lemma 4.1}
Suppose that $A_1,A_2,\ldots, A_n$ are some random events under an
identical probability space such that $P(A_i)>0$ for $1\leq i\leq
n$, then
\[
P\big(\bigcup\limits_{i=1}^nA_i\big)\geq\frac{1}{\frac{1}{n^2}\sum\limits_{i=1}^n
\sum\limits_{j=1}^n\frac{P(A_i\bigcap A_j)}{P(A_i)P(A_j)}}.
\]
\end{lemma}
\proof For $1\leq i\leq n$, we define $Y_i=\frac{1_{A_i}}{P(A_i)}$,
then
\[
\{\sum\limits_{i=1}^nY_i>0\}=\bigcup\limits_{i=1}^nA_i.
\]
As a result, according to H\"{o}lder's inequality,
\begin{align}\label{equ 4.3}
P(\bigcup\limits_{i=1}^nA_i)&=P(\sum\limits_{i=1}^nY_i>0)\geq\frac{[{\rm E}\sum\limits_{i=1}^nY_i]^2}{{\rm E}(\sum\limits_{i=1}^nY_i)^2}\notag\\
&=\frac{n^2}{\sum\limits_{i=1}^n\sum\limits_{j=1}^n{\rm
E}(Y_iY_j)}=\frac{n^2}{\sum\limits_{i=1}^n\sum\limits_{j=1}^n\frac{P(A_i\bigcap
A_j)}{P(A_i)P(A_j)}}.
\end{align}

\qed

To give a crucial lemma for the proof of Equation \eqref{equ 2.2},
we introduce following definitions about a random walk on
$\mathbb{Z}^d$. We define $\{S_n\}_{n\geq 0}$ as a random walk on
$\mathbb{Z}^d$ such that
\[
P(S_{j+1}-S_j=e_k)=\frac{1}{\lfloor \frac{d}{N}\rfloor}
\]
for $d-\lfloor \frac{d}{N}\rfloor+1\leq k\leq d$ and $j$ satisfying
that $N|(j+1)$ while
\[
P(S_{j+1}-S_j=e_l)=P(S_{j+1}-S_j=-e_l)=\frac{1}{2(d-\lfloor
\frac{d}{N}\rfloor)}
\]
for $1\leq l\leq d-\lfloor \frac{d}{N}\rfloor$ and $j$ satisfying
that $N\not|(j+1)$. This random walk is first introduced in
\cite{Kesten1990} by Kesten. We denote by $\{\widehat{S}_n\}_{n\geq
0}$ an independent copy of $\{S_n\}_{n\geq 0}$. For any $x,y\in
\mathbb{Z}^d$, we denote by $\widetilde{P}_{x,y}$ the probability
measure of $\{S_n,\widehat{S}_n\}_{n\geq 0}$ with
$S_0=x,\widehat{S}_0=y$. We denote by $\widetilde{\rm E}_{x,y}$ the
expectation operator with respect to $\widetilde{P}_{x,y}$. For each
$x\in \mathbb{Z}^d$, we define
\[
K(x,S)=|\{i\geq 0:S_i=x\}|
\]
as the times that $S$ visits $x$. Similarly, we define
\[
K(x,\widehat{S})=|\{i\geq 0:\widehat{S}_i=x\}|.
\]
Let $L(x,S,\widehat{S})=\min\{K(x,S),K(x,\widehat{S})\}$ and
\[
L(S,\widehat{S})=\sum_{x\in \mathbb{Z}^d}L(x,S,\widehat{S}),
\]
then the following lemma is crucial for us to prove Equation
\eqref{equ 2.2}.
\begin{lemma}\label{lemma 4.2}
\begin{equation}\label{equ 4.4}
P_{\lambda,d}(C_t^{A}\neq \emptyset,\forall~t>0)\geq
\frac{1}{\frac{1}{|A|^2}\sum\limits_{x\in A}\sum\limits_{y\in
A}\widetilde{\rm E}_{x,y}\big(\frac{1}{q}\big)^{L(S,\widehat{S})}},
\end{equation}
where
\[
q=\frac{\lambda p}{2d}{\rm
E}\big(\frac{1}{\xi+\frac{\lambda}{2d}}\big).
\]
\end{lemma}

We will give the proof of Lemma \ref{lemma 4.2} later. First we show
that how to utilize Lemma \ref{lemma 4.2} to prove Equation
\eqref{equ 2.2}.

\proof[Proof of Equation \eqref{equ 2.2}]

For $\lambda=\gamma\lambda_c$ with $\gamma>1$,
\begin{equation}\label{equ 4.4 two}
2dq\geq \widetilde{\gamma}=\frac{\gamma+1}{2}
\end{equation}
for sufficiently large $d$. Let $H$ be the event that there exists
$x\in \mathbb{Z}^d$ such that $L(x,S,\widehat{S})\geq 1$, then
\begin{equation}\label{equ 4.5}
\widetilde{\rm E}_{x,y}\big(\frac{1}{q}\big)^{L(S,\widehat{S})}
=\widetilde{P}_{x,y}(H^c)+\widetilde{\rm E}_{x,y}1_H\Big({\rm
E}\big((\frac{1}{q})^{L(S,\widehat{S})}|H\big)\Big).
\end{equation}
In \cite{Kesten1990}, Kesten gives a detailed calculation of the
upper bound of generating function of $L(S,\widehat{S})$ (which is
denoted by $J(r,r^\prime)$ in that paper). Due to the analysis in
\cite{Kesten1990} which leads to Equation (2.44) and Lemma 7 of that
paper,
\begin{equation}\label{equ 4.6}
{\rm E}\big((\frac{1}{q})^{L(S,\widehat{S})}|H\big)\leq
M_6\sum_{n=0}^{+\infty}\beta(d)^n,
\end{equation}
where
\[
\beta(d)=\frac{N^{\frac{3}{N-1}}}{2q(d-\lfloor\frac{d}{N}\rfloor)}+\frac{12}{qdN^2}
+3M_7\frac{N^5}{qd^2}+3M_7q^{-N}(\frac{N}{2(d-\lfloor\frac{d}{N}\rfloor)})^{N+1}N^3
\]
and $M_6, M_7$ are constants which do not depend on $d$.

According the definition of $\beta(d)$ and Equation \eqref{equ 4.4
two},
\begin{equation}\label{equ 4.8}
\limsup_{d\rightarrow+\infty}\beta(d)\leq
\frac{1}{\widetilde{\gamma}}<1.
\end{equation}
By Equation \eqref{equ 4.6} and \eqref{equ 4.8}, for sufficiently
large $d$,
\begin{equation}\label{equation 4.9}
{\rm E}\big((\frac{1}{q})^{L(S,\widehat{S})}|H\big)<M_8<+\infty,
\end{equation}
where $M_8$ is a constant which does not depend on $d$.

By Equation \eqref{equ 4.5} and \eqref{equation 4.9}, for
$\lambda>\lambda_c$ and sufficiently large $d$,
\begin{equation}\label{equ 4.10}
\widetilde{\rm E}_{x,y}\big(\frac{1}{q}\big)^{L(S,\widehat{S})}\leq
1-\widetilde{P}_{x,y}(H)+\widetilde{P}_{x,y}(H)M_8\leq
1+\widetilde{P}_{x,y}(H)M_8.
\end{equation}
When $x=y$, $\widetilde{P}_{x,y}(H)=1$ since
$\{S_0=\widehat{S}_0\}$. When $x\neq y$, we claim that
\begin{equation}\label{equ 4.10 two}
\widetilde{P}_{x,y}(H)\leq \frac{M_9 N}{d},
\end{equation}
where $M_9$ is a constant which does not depend on $d,x,y$. If
Equation \eqref{equ 4.10 two} holds, then by Lemma \ref{lemma 4.2},
Equation \eqref{equ 4.10} and \eqref{equ 4.10 two},
\begin{align}\label{equ 4.13}
P_{\lambda,d}(C_t^{A(d)}\neq \emptyset,\forall~t>0)&\geq
\frac{1}{\frac{1}{\log^2 d}[(\log d)M_8+(\log^2 d-\log
d)(1+\frac{M_8M_9N}{d})]}\notag\\
&=\frac{1}{\frac{M_8}{\log d}+(1-\frac{1}{\log
d})(1+\frac{M_8M_9N}{d})}.
\end{align}
Equation \eqref{equ 2.2} follows from \eqref{equ 4.13} directly.

To finish the proof, we only need to show that Equation \eqref{equ
4.10 two} holds. For $x\neq y$,
\begin{equation}\label{equ 4.11}
H=\{\exists~ i,j\geq 0, i+j>0, S_i=\widehat{S}_j\}.
\end{equation}
Therefore,
\begin{equation}\label{equ 4.12}
\widetilde{P}_{x,y}(H)\leq
\sum_{i=1}^{+\infty}\widetilde{P}_{x,y}(\exists~j\geq 0,
S_i=\widehat{S}_j)+\sum_{j=1}^{+\infty}\widetilde{P}_{x,y}(\exists~i\geq
0, \widehat{S}_j=S_i).
\end{equation}
For any $x=(x_1,x_2,\ldots,x_d)\in \mathbb{Z}^d$, we define
\[
\sigma(x)=\sum_{j=d-\lfloor\frac{d}{N}\rfloor+1}^dx_j,
\]
then according to the definition of $S$ and $\widehat{S}$,
\[
\sigma(S_n)=\sigma(S_0)+\lfloor\frac{n}{N}\rfloor,\text{\quad}
\sigma(\widehat{S}_n)=\sigma(\widehat{S}_0)+\lfloor\frac{n}{N}\rfloor
\]
for any $n\geq 1$. As a result, for each $i\geq 1$, if there exists
$j$ such that $S_i=\widehat{S}_j$, then
\begin{equation}\label{equ 4.15}
N\Big(\sigma(S_0)-\sigma(\widehat{S}_0)+\lfloor\frac{i}{N}\rfloor\Big)\leq
j\leq
N\Big(\sigma(S_0)-\sigma(\widehat{S}_0)+\lfloor\frac{i}{N}\rfloor\Big)+N-1.
\end{equation}
As a result,
\begin{equation}\label{equ 4.16}
\widetilde{P}_{x,y}(\exists~j\geq 0, S_i=\widehat{S}_j) \leq
\sum_{j=N\Big(x-y+\lfloor\frac{i}{N}\rfloor\Big)}^{N\Big(x-y+\lfloor\frac{i}{N}\rfloor\Big)+N-1}\widetilde{P}_{x,y}(S_i=\widehat{S}_j)
\leq N\sup_{u\in \mathbb{Z}^d}\widetilde{P}_{x}(S_i=u).
\end{equation}
By Equation \eqref{equ 4.16},
\begin{align}\label{equ 4.17}
\sum_{i=1}^{+\infty}\widetilde{P}_{x,y}(\exists~j\geq 0,
S_i=\widehat{S}_j)&\leq N\sum_{i=1}^{+\infty}\sup_{u\in
\mathbb{Z}^d}\widetilde{P}_{x}(S_i=u)\notag\\
&=N\sup_{u\in
\mathbb{Z}^d}\widetilde{P}_{x}(S_1=u)+N\sum_{i=2}^{+\infty}\sup_{u\in
\mathbb{Z}^d}\widetilde{P}_{x}(S_i=u)\notag\\
&=\frac{N}{2(d-\lfloor\frac{d}{N}\rfloor+1)}+N\sum_{i=2}^{+\infty}\sup_{u\in
\mathbb{Z}^d}\widetilde{P}_{x}(S_i=u).
\end{align}
For $x=(x_1,x_2,\ldots,x_d)\in \mathbb{Z}^d$, we define
\[
\beta(x)=(x_1,x_2,\ldots,x_{d-\lfloor\frac{d}{N}\rfloor+1})\in
\mathbb{Z}^{d-\lfloor\frac{d}{N}\rfloor+1}.
\]
We define $\{\phi_n\}_{n\geq 0}$ as the simple random walk on
$\mathbb{Z}^{d-\lfloor\frac{d}{N}\rfloor+1}$. According to the
definition of $S$, $\beta(S_n)$ and
$\phi_{n-\lfloor\frac{n}{N}\rfloor}$ have identical probability
distribution when $\beta(S_0)=\phi_0$. As a result,
\begin{align}\label{equ 4.18}
\sum_{i=2}^{+\infty}\sup_{u\in \mathbb{Z}^d}\widetilde{P}_{x}(S_i=u)
&\leq \sum_{i=2}^{+\infty}\sup_{u\in
\mathbb{Z}^{d-\lfloor\frac{d}{N}\rfloor+1}}\widetilde{P}_{\beta(x)}(\phi_{i-\lfloor\frac{i}{N}\rfloor}=u)\notag\\
&\leq 2\sum_{i=2}^{+\infty}\sup_{u\in
\mathbb{Z}^{d-\lfloor\frac{d}{N}\rfloor+1}}\widetilde{P}_{\beta(x)}(\phi_i=u).
\end{align}
For simple random walk $\{X_n\}_{n\geq 0}$ on $\mathbb{Z}^d$, it is
shown in \cite{Kesten1964} that
\begin{equation}\label{equ 4.19}
\sum_{i=2}^{+\infty}\sup_{u\in \mathbb{Z}^d}P_x(X_i=u) \leq
2\sum_{i=1}^{+\infty}P_x(X_{2i}=x)\leq \frac{M_{10}}{d},
\end{equation}
where $M_{10}$ is a constant which does not depend on $d$. By
Equation \eqref{equ 4.19},
\begin{equation}\label{equ 4.20}
\sum_{i=2}^{+\infty}\sup_{u\in
\mathbb{Z}^{d-\lfloor\frac{d}{N}\rfloor+1}}\widetilde{P}_{\beta(x)}(\phi_i=u)
\leq \frac{M_{10}}{d-\lfloor\frac{d}{N}\rfloor+1}.
 \end{equation}
By Equation \eqref{equ 4.17}, \eqref{equ 4.18} and \eqref{equ 4.20},
\begin{equation}\label{equ 4.21}
\sum_{i=1}^{+\infty}\widetilde{P}_{x,y}(\exists~j\geq 0,
S_i=\widehat{S}_j)\leq \frac{NM_{11}}{d}
\end{equation}
for sufficiently large $d$ and $x\neq y$, where $M_{11}$ is a
constant which does not depend on $d,x,y$.

Equation \eqref{equ 4.10 two} follows from \eqref{equ 4.12} and
\eqref{equ 4.21} directly and the proof of Equation \eqref{equ 2.2}
is complete.

\qed

Now we give the proof of Lemma \ref{lemma 4.2}.

\proof[Proof of Lemma \ref{lemma 4.2}]

For each $k\geq 1$ and $x\in A$,
\[
|F_k(x)|=g_k=2^{(N-1)k}(d-\lfloor\frac{d}{N}\rfloor)^{(N-1)k}{\lfloor\frac{d}{N}\rfloor}^{k-1}.
\]

Then, by Lemma \ref{lemma 4.1},
\begin{equation}\label{equ 4.22}
P_{\lambda,d}\Big(\bigcup\limits_{x\in A}\bigcup\limits_{\vec{x}\in
F_k(x)}I_k(\vec{x})\Big)\geq
\frac{1}{\frac{1}{g_k^2|A|^2}\sum\limits_{x,y\in
A}\sum\limits_{\vec{x}\in F_k(x),\atop\vec{y}\in
F_k(y)}\frac{P_{\lambda,d}(I_k(\vec{x})\cap
I_k(\vec{y}))}{P_{\lambda,d}(I_k(\vec{x}))P_{\lambda,d}(I_k(\vec{y}))}}.
\end{equation}
According to the definition of $I_k(\vec{x})$,
\begin{align}\label{equ 4.23}
P_{\lambda,d}(I_k(\vec{x}))&={\rm
E}_{\lambda,d}[P_\lambda^\omega(W_i\leq
s_i,0\leq i\leq kN-2)\prod\limits_{i=0}^{kN-2}X(x_i,x_{i+1})]\notag\\
&={\rm
E}_{\lambda,d}\Big(\prod\limits_{i=0}^{kN-2}X(x_i,x_{i+1}){\rm
E}_\lambda^\omega(\frac{\frac{\lambda}{2d}}{\frac{\lambda}{2d}+\xi(x_i)})\Big),
\end{align}
where
\[
\vec{x}=(x_0,x_1,\ldots,x_{kN-1})
\]
while $\{W_i\}_{0\leq i\leq kN-1}$ are i. i. d. exponential times
with rate $\frac{\lambda}{2d}$ and $\{s_i\}_{0\leq i\leq kN-1}$ are
random exponential times with rate
$\xi(x_0),\xi(x_1),\ldots,\xi(x_{kN-2})$ respectively. Please note
that the factor $\prod\limits_{i=0}^{kN-2}X(x_i,x_{i+1})$ in
Equation \eqref{equ 4.23} is the index of the event that all the
edges on the path $\vec{x}$ are open.

By Equation \eqref{equ 4.23},
\begin{align}\label{equ 4.24}
P_{\lambda,d}(I_k(\vec{x}))P_{\lambda,d}(I_k(\vec{y}))=&{\rm
E}_{\lambda,d}\Big(\prod\limits_{i=0}^{kN-2}X(x_i,x_{i+1})\Big){\rm
E}_{\lambda,d}\Big(\prod\limits_{i=0}^{kN-2}X(y_i,y_{i+1})\Big)\\
&\times{\rm E}_{\lambda,d}\Big(\prod\limits_{i=0}^{kN-2} {\rm
E}_\lambda^\omega(\frac{\frac{\lambda}{2d}}{\frac{\lambda}{2d}+\xi(x_i)})\Big){\rm
E}_{\lambda,d}\Big(\prod\limits_{i=0}^{kN-2}{\rm
E}_\lambda^\omega(\frac{\frac{\lambda}{2d}}{\frac{\lambda}{2d}+\xi(y_i)})\Big)\notag,
\end{align}
since $\{\xi(x)\}_{x\in \mathbb{Z}^d}$ and $\{X(e)\}_{e\in
\mathbb{E}^d}$ are independent as we assumed.

For any $u\in \mathbb{Z}^d,\vec{x}\in F_k(x)$ and $\vec{y}\in
F_k(y)$, we define
\begin{equation}\label{equ 4.25}
l_k(u,\vec{x})=\{0\leq i\leq kN-2:x_i=u\}.
\end{equation}
and
\begin{equation}\label{equ 4.26}
G(u,\vec{x},\vec{y})=
\begin{cases}
\Big[{\rm
E}_\lambda^{\omega}(\frac{\frac{\lambda}{2d}}{\frac{\lambda}{2d}+\xi(u)})\Big]^{l_k(u,\vec{x})}&
\text{\quad if~}|l_k(u,\vec{x})|\geq |l_k(u,\vec{y})|,\\
\Big[{\rm
E}_\lambda^{\omega}(\frac{\frac{\lambda}{2d}}{\frac{\lambda}{2d}+\xi(u)})\Big]^{l_k(u,\vec{y})}&
\text{\quad if~}|l_k(u,\vec{x})|< |l_k(u,\vec{y})|.
\end{cases}
\end{equation}
Then, it is not difficult to see that
\begin{equation}\label{equ 4.27}
P_{\lambda}^\omega(I_k(\vec{x})\cap I_k(\vec{y}))\leq
\Big[\prod\limits_{i=0}^{kN-2}X(x_i,x_{i+1})X(y_i,y_{i+1})\Big]\prod_{u\in
\mathbb{Z}^d}G(u,\vec{x},\vec{y}).
\end{equation}
The explanation of Equation \eqref{equ 4.27} is that when the times
$\vec{x}$ visits $u$ is bigger than that of $\vec{y}$, then we do
not care the probability that $y_j$ infects $y_{j+1}$ for each $j$
such that $y_j=u$, then we will obtain an upper bound of
$P_{\lambda}^\omega(I_k(\vec{x})\cap I_k(\vec{y}))$.

We define
\begin{equation*}
m(\vec{x},\vec{y})=\{u\in \mathbb{Z}^d: |l_k(u,\vec{x})|\geq
|l_k(u,\vec{y})|\}
\end{equation*}
and
\begin{equation*}
n(\vec{x},\vec{y})=\{u\in \mathbb{Z}^d: |l_k(u,\vec{x})|<
|l_k(u,\vec{y})|\},
\end{equation*}
then by Equation \eqref{equ 4.27},
\begin{align}\label{equ 4.28}
&P_{\lambda,d}(I_k(\vec{x})\cap I_k(\vec{y}))\notag\\
&\leq {\rm E}_{\lambda,d}\big[\prod\limits_{i=0}^{kN-2}X(x_i,x_{i+1})X(y_i,y_{i+1})\big]\notag\\
&\times {\rm E}_{\lambda,d}\Big(\prod\limits_{u\in
m(\vec{x},\vec{y})}\prod\limits_{j\in l_k(u,\vec{x})}{\rm
E}_\lambda^{\omega}(\frac{\frac{\lambda}{2d}}{\frac{\lambda}{2d}+\xi(u)})\Big){\rm
E}_{\lambda,d}\Big(\prod\limits_{u\in
n(\vec{x},\vec{y})}\prod\limits_{j\in l_k(u,\vec{y})}{\rm
E}_\lambda^{\omega}(\frac{\frac{\lambda}{2d}}{\frac{\lambda}{2d}+\xi(u)})\Big),
\end{align}
since $\{\xi(u):u\in m(\vec{x},\vec{y})\}$ and $\{\xi(v):v\in
n(\vec{x},\vec{y})\}$ are independent.

According to Equation \eqref{equ 4.24} and \eqref{equ 4.28},
\begin{equation}\label{equ 4.29}
\frac{P_{\lambda,d}(I_k(\vec{x})\cap
I_k(\vec{y}))}{P_{\lambda,d}(I_k(\vec{x}))P_{\lambda,d}(I_k(\vec{y}))}\leq
\Gamma_1\Gamma_2,
\end{equation}
where
\[
\Gamma_1=\frac{{\rm
E}_{\lambda,d}\big(\prod\limits_{i=0}^{kN-2}X(x_i,x_{i+1})X(y_i,y_{i+1})\big)}{{\rm
E}_{\lambda,d}\Big(\prod\limits_{i=0}^{kN-2}X(x_i,x_{i+1})\Big){\rm
E}_{\lambda,d}\Big(\prod\limits_{i=0}^{kN-2}X(y_i,y_{i+1})\Big)}
\]
and
\[
\Gamma_2=\frac{1}{{\rm E}_{\lambda,d}\Big(\prod\limits_{u\in
n(\vec{x},\vec{y})}\prod\limits_{j\in l_k(u,\vec{x})}{\rm
E}_\lambda^{\omega}(\frac{\frac{\lambda}{2d}}{\frac{\lambda}{2d}+\xi(u)})\Big){\rm
E}_{\lambda,d}\Big(\prod\limits_{u\in
m(\vec{x},\vec{y})}\prod\limits_{j\in l_k(u,\vec{y})}{\rm
E}_\lambda^{\omega}(\frac{\frac{\lambda}{2d}}{\frac{\lambda}{2d}+\xi(u)})\Big)}
\]
since $m(\vec{x},\vec{y})\cup n(\vec{x},\vec{y})=\mathbb{Z}^d$,
$\{\xi(u):u\in m(\vec{x},\vec{y})\}$ and $\{\xi(v):v\in
n(\vec{x},\vec{y})\}$ are independent and
\[
{\rm E}_{\lambda,d}\Big(\prod\limits_{i=0}^{kN-2} {\rm
E}_\lambda^\omega(\frac{\frac{\lambda}{2d}}{\frac{\lambda}{2d}+\xi(x_i)})\Big)=
{\rm E}_{\lambda,d}\Big(\prod\limits_{u\in
\mathbb{Z}^d}\prod\limits_{j\in l_k(u,\vec{x})}{\rm
E}_\lambda^{\omega}(\frac{\frac{\lambda}{2d}}{\frac{\lambda}{2d}+\xi(u)})\Big)
\]
Since $\{\xi(x)\}_{x\in \mathbb{Z}^d}$ are i. i. d., $\{{\rm
E}_\lambda^{\omega}(\frac{\frac{\lambda}{2d}}{\frac{\lambda}{2d}+\xi(x_j)})\}_{u\in
n(\vec{x},\vec{y}),j\in l_k(u,\vec{x})}$ are positive correlated. As
a result,
\begin{align*}
{\rm E}_{\lambda,d}\Big(\prod\limits_{u\in
n(\vec{x},\vec{y})}\prod\limits_{j\in l_k(u,\vec{x})}{\rm
E}_\lambda^{\omega}(\frac{\frac{\lambda}{2d}}{\frac{\lambda}{2d}+\xi(u)})\Big)
&\geq \Big({\rm
E}_{\lambda,d}(\frac{\frac{\lambda}{2d}}{\frac{\lambda}{2d}+\xi(u)})\Big)^{\sum\limits_{u\in
n(\vec{x},\vec{y})}|l_k(u,\vec{x})|}\\
&=\big(\frac{q}{p}\big)^{\sum\limits_{u\in
n(\vec{x},\vec{y})}|l_k(u,\vec{x})|}
\end{align*}
and
\[
{\rm E}_{\lambda,d}\Big(\prod\limits_{u\in
m(\vec{x},\vec{y})}\prod\limits_{j\in l_k(u,\vec{y})}{\rm
E}_\lambda^{\omega}(\frac{\frac{\lambda}{2d}}{\frac{\lambda}{2d}+\xi(u)})\Big)\geq
\big(\frac{q}{p}\big)^{\sum\limits_{u\in
m(\vec{x},\vec{y})}|l_k(u,\vec{y})|},
\]
where $q$ is defined in \eqref{equ 4.4}.

Therefore,
\begin{equation}\label{equ 4.30}
\Gamma_2\leq \big(\frac{p}{q}\big)^{L_k(\vec{x},\vec{y})},
\end{equation}
where
\[
L_k(\vec{x},\vec{y})=\sum\limits_{u\in
n(\vec{x},\vec{y})}|l_k(u,\vec{x})|+\sum\limits_{u\in
m(\vec{x},\vec{y})}|l_k(u,\vec{y})|=\sum_{u\in
\mathbb{Z}^d}\min\{|l_k(u,\vec{x})|,|l_k(u,\vec{y})|\}.
\]

We define
\[
\Lambda=\{e\in \mathbb{E}^d: \text{both $\vec{x}$ and $\vec{y}$
visit $e$}\},
\]
then it is easy to see that
\begin{equation}\label{equ 4.31}
\Gamma_1=\big(\frac{1}{p}\big)^{|\Lambda|}.
\end{equation}
For each $x\in \mathbb{Z}^d$, we define
\[
\Lambda(x)=\{e\in \Lambda:e\text{~connects~}x\},
\]
then
\[
2|\Lambda|=\sum_{x\in \mathbb{Z}^d}|\Lambda(x)|.
\]
Since in a path each vertex connects two edges,
\begin{equation}\label{equ 4.32}
\frac{|\Lambda(u)|}{2}\leq \min\{|l_k(u,\vec{x})|,|l_k(u,\vec{y})|\}
\end{equation}
for each $u\in \mathbb{Z}^d$. Therefore,
\begin{equation}\label{equ 4.33}
|\Lambda|=\sum_{u\in \mathbb{Z}^d}\frac{|\Lambda(u)|}{2}\leq
\sum_{u\in
\mathbb{Z}^d}\min\{|l_k(u,\vec{x})|,|l_k(u,\vec{y})|\}=L_k(\vec{x},\vec{y}).
\end{equation}
By Equation \eqref{equ 4.31} and \eqref{equ 4.33},
\begin{equation}\label{equ 4.34}
\Gamma_1\leq \big(\frac{1}{p}\big)^{L_k(\vec{x},\vec{y})}.
\end{equation}
By Equation \eqref{equ 4.29}, \eqref{equ 4.30} and \eqref{equ 4.34},
\begin{equation}\label{equ 4.35}
\frac{P_{\lambda,d}(I_k(\vec{x})\cap
I_k(\vec{y}))}{P_{\lambda,d}(I_k(\vec{x}))P_{\lambda,d}(I_k(\vec{y}))}\leq
\big(\frac{1}{q}\big)^{|L_k(\vec{x},\vec{y})|}.
\end{equation}
By Equation \eqref{equ 4.22} and \eqref{equ 4.35},
\begin{align}\label{equ 4.36}
P_{\lambda,d}\Big(\bigcup\limits_{x\in A}\bigcup\limits_{\vec{x}\in
F_k(x)}I_k(\vec{x})\Big)&\geq
\frac{1}{\frac{1}{g_k^2|A|^2}\sum\limits_{x,y\in
A}\sum\limits_{\vec{x}\in F_k(x),\atop\vec{y}\in
F_k(y)}\big(\frac{1}{q}\big)^{|L_k(\vec{x},\vec{y})|}}\\
&=\frac{1}{\frac{1}{|A|^2}\sum\limits_{x,y\in A}\widetilde{\rm
E}_{x,y}\big(\frac{1}{q}\big)^{|L_k(S,\widehat{S})|}},\notag
\end{align}
where
\[
L_k(S,\widehat{S})=L_k(\{S_i\}_{i\leq kN-1},\{\widehat{S}_i\}_{i\leq
kN-1}).
\]
According to the definition of $L(S,\widehat{S})$ before Lemma
\eqref{lemma 4.2},
\begin{equation}\label{equ 4.37}
\lim_{k\rightarrow+\infty}L_k(S,\widehat{S})=L(S,\widehat{S}).
\end{equation}
Lemma \ref{lemma 4.2} follows from Equation \eqref{equ 4.2}
,\eqref{equ 4.36} and \eqref{equ 4.37} directly.

\qed

At last we give the proof of Corollary \ref{corollary 2.2}.

\proof[Proof of Corollary \ref{corollary 2.2}]

By Equation \eqref{equ 2.2}, for any $\lambda>\lambda_c$,
\begin{equation}\label{equ 4.38}
P_{\lambda,d}(C_t^{A(d)}\neq \emptyset,\forall~t>0)>0
\end{equation}
for sufficiently large $d$. By Equation \eqref{equ 3.3},
\begin{align}\label{equ 4.39}
P_{\lambda,d}(C_t^{O}\neq
\emptyset,\forall~t>0)&=\lim\limits_{t\rightarrow+\infty}P_{\lambda,d}(C_t^{O}\neq
\emptyset)\notag\\
&\geq
\frac{1}{|A(d)|}\lim\limits_{t\rightarrow+\infty}P_{\lambda,d}(C_t^{A(d)}\neq
\emptyset)=\frac{P_{\lambda,d}(C_t^{A(d)}\neq
\emptyset,\forall~t>0)}{|A(d)|}.
\end{align}

By Equation \eqref{equ 4.38} and \eqref{equ 4.39},
\begin{equation}\label{equ 4.40}
P_{\lambda,d}(C_t^{O}\neq \emptyset,\forall~t>0)>0
\end{equation}
for sufficiently large $d$. By the definition of $\lambda_d$ in
Equation \eqref{equ 2.3},
\[
\lambda_d\leq \lambda
\]
for sufficiently large $d$ and hence
\[
\limsup_{d\rightarrow+\infty}\lambda_d\leq \lambda.
\]
Let $\lambda\rightarrow \lambda_c$, then the proof is complete.

\qed

\quad

\textbf{Acknowledgments.} The author is grateful to the financial
support from the National Natural Science Foundation of China with
grant number 11501542 and China Postdoctoral Science Foundation (No.
2015M571095).

{}
\end{document}